\documentclass[12pt,reqno]{article}
\def\hybrid{\topmargin      0pt
\oddsidemargin 0pt
\headheight 0pt \headsep 0pt
\textwidth 160true mm
\textheight 231true mm
\marginparwidth 0.0in
\parskip 3pt plus 1pt   \jot = 1.5ex}
\usepackage{amssymb}
\usepackage{amsmath}
\usepackage{amscd}
\usepackage{amsthm}
\hybrid

\theoremstyle{plain} 
\newtheorem{thm}{Theorem}[section]
\newtheorem{cor}[thm]{Corollary}
\newtheorem{lemma}[thm]{Lemma}
\newtheorem{propn}[thm]{Proposition}

\theoremstyle{definition}

\theoremstyle{definition}  

\newcommand{\C}{\mathbb{C}}

\newcommand{\Z}{\mathbb{Z}}

\newcommand{\ff}{\varphi}

\newcommand{\tE}{{\widetilde E}}
\renewcommand{\P}{\mathbb{P}}

\newcommand{\LL}{\mathcal{L}}

\newcommand{\OO}{\mathcal{O}}
\renewcommand{\O}{\mathcal{O}}

\newcommand{\End}{\mbox{End}}

\renewcommand{\S}{\mathcal{S}}

\newcommand{\be}[1]{\begin{eqnarray#1}}
\newcommand{\ee}[1]{\end{eqnarray#1}}
\newcommand{\tor}[1]{\stackrel{#1}{\longrightarrow}}
\newcommand{\ot}{\otimes}
\newcommand{\Pic}{\mbox{Pic}}
\newcommand{\codim}{\mbox{codim}}

\renewcommand{\span}{{\rm span}}
\DeclareMathOperator{\rk}{rk}

\begin{document}

\title{Finite rank vector bundles on inductive limits of grassmannians}
\author{Joseph Donin\thanks{Work supported in part
by Israel Academy of Sciences Grant no. 8007/99-03 and the Max Planck Institute
for Mathematics, Bonn}
\ and Ivan Penkov\thanks{Work supported in part
by an NSF grant and the Max Planck Institute
for Mathematics, Bonn}}

\date{}
\maketitle
\begin{abstract}
If $\P^\infty$ is the projective ind-space, i.e. $\P^\infty$ is the
inductive limit  of linear embeddings of complex
projective spaces,
the Barth-Van de Ven-Tyurin (BVT)
 Theorem claims that every finite rank vector bundle on $\P^\infty$
is isomorphic to a direct sum of line bundles. We extend this
theorem to general sequences of morphisms between projective
spaces by proving that, if there are infinitely many morphisms of
degree higher than one, every vector bundle of finite rank on the
inductive limit is trivial. We then establish a relative version
of these results, and apply it to  the study of vector bundles on
inductive limits of grassmannians. In particular we show that the
BVT Theorem extends to the ind-grassmannian of subspaces
commensurable with a fixed infinite dimensional and infinite
codimensional subspace in $\C^\infty$. We also show that for a
class of twisted ind-grassmannians, every finite rank vector
bundle is trivial.

2000 AMS Subject Classification: Primary 32L05, 14J60,  Secondary 14M15.
\end{abstract}

\vskip-5.6in

\hskip3.75in
{\small To the memory of Andrey Tyurin}

\vskip5.6in

\section{Introduction}
About 30 years ago the study of ``infinitely extendable''
vector bundles of finite rank
on projective spaces and grassmannians was initiated.
In particular the following remarkable
theorem was proved: any finite rank vector bundle on the
infinite complex projective space $\P^\infty$ (or equivalently,
any finite rank vector bundle on $\P^n$ which admits an extension
to $\P^m$ for large enough $m>n$) is isomorphic to a direct sum
of line bundles. For rank two bundles this was established by
W. Barth and A. Van de Ven in \cite{BV}, and for arbitrary finite
rank bundles the theorem was proved by A. Tyurin in \cite{T}.
In what follows we refer to this result as to the
Barth-Van de Ven-Tyurin Theorem, or as to the BVT Theorem.
Some first steps were made also towards understanding finite rank
vector bundles on infinite grassmannians:
R. Hartshorne conjectured that every finite rank bundle on an infinite
grassmannian $G(k,\infty)$ is homogeneous, see \cite{BV}.
For rank two bundles this conjecture
is proved  in \cite{BV}, and in the general case the conjecture is proved by E. Sato, \cite{S2}.
Sato established also a partial analog of Hartshorne's conjecture for the infinite grassmanians of the classical simple groups, and reproved the BVT Theorem, see \cite{S2} and \cite{S1}.

The purpose of the present note is to revive this discussion
and relate it with the more recent duscussion of homogeneous ind-spaces of
locally linear ind-groups, see \cite{DPW}, \cite{DP} and the references therein.
Our starting point is an infinite sequence
\be{*}
G_1\subset G_2\subset...
\ee{*}
of complex linear algebraic groups and a subsequence
\be{*}
P_1\subset P_2\subset...
\ee{*}
of parabolic subgroups. This yields a sequence of morphisms
\be{}\label{eq1}
G_1/P_1\tor{} G_2/P_2\tor{}...\ .
\ee{}
Let $G/P$ denote the inductive limit of (\ref{eq1}).
We restrict ourselves here to the case when $G_N\simeq GL(n_N)$
and $P_N$ are maximal parabolic subgroups, i.e. $G_N/P_N$ are grassmannians.

The study of line
bundles on $G/P$ gives some first hints on what
general finite rank bundles on $G/P$ might look like.
An essential difference with the cases studied in \cite{BV}, \cite{T}, \cite{S1}, and \cite{S2}
is that the restriction maps
\be{*}
\Pic(G_N/P_N)\tor{}\Pic(G_{N-1}/P_{N-1})
\ee{*}
on the Picard groups induced by inclusions in (\ref{eq1}) are only
injective and in general not surjective.
Therefore $\Pic(G/P)$, the Picard group of the inductive limit,
is isomorphic to $\Z$ or equals zero. In the first case the call the sequence
(\ref{eq1}) {\em linear} and in the second case we call it {\em twisted}.

Here is a brief description of our results.
The simplest case for a sequence (\ref{eq1}) is when all $G_N/P_N$
are projective spaces. Consider more generally an arbitrary sequence of morphisms of
projective spaces
\be{}\label{eq2}
\P^{i_1}\tor{\ff_1} \P^{i_2}\tor{\ff_2} ...\ .
\ee{}
If (\ref{eq2}) is linear, i.e.  all but finitely many Picard groups map isomorphically,
the inductive limit is isomorphic to the projective ind-space
$\P^\infty$.
Here the BVT Theorem claims that any finite rank bundle is isomorphic
to a direct sum of line bundles.
If (\ref{eq2}) is twisted, i.e. the Picard group of the inductive limit $\P^\infty_{tw}$
of (\ref{eq2}) equals zero,
the problem of describing all finite rank vector bundles on $\P^\infty_{tw}$
was not posed in the 70's.

Our first main result (Theorem \ref{thmour}) claims that every
finite rank vector bundle on $\P^\infty_{tw}$ is trivial.
This theorem, together with the BVT Theorem, can be extended to the relative
case  and yields
a complete description of finite rank vector bundles on any inductive
limit of relative projective spaces.

In Section 4 we apply the above result to the study of finite rank vector bundles
on inductive limits of grassmannians $G(k_N,n_N)$. We consider two types of morphisms
of grassmannians in (\ref{eq1}): standard inclusions and a certain class of twisted
homogeneous morphisms which we call twisted extensions.
For standard inclusions we show that, if
$\lim_{N\to\infty}k_N=\lim_{N\to\infty}(n_N-k_N)=\infty$, any finite rank
vector bundle on the inductive limit is isomorphic to a direct sum
of line bundles. An interesting ind-variety which arises as
the inductive limit of a sequence of standard inclusions
satisfying the above condition is the ind-grassmanian
$G(V, \infty)$ of subspaces $V' \subset \C^{\infty}$
commensurable with a fixed infinite dimensional and infinite
codimensional subspace
$V \subset \C^{\infty}$. Therefore, any finite rank vector
bundle on $G(V, \infty)$ is isomorphic to a direct sum of line bundles.
Finally, we prove that in the case of twisted extensions of grassmannians every
finite rank vector bundle on the inductive limit is trivial.

{\bf Acknowledgements.} We are grateful to A. Tyurin for a helpful discussion on
the topic of this paper several weeks before he suddenly passed away in October 2002.
We thank also Z. Ran for his comments and for making us aware of E. Sato's work.

\section{Preliminary results}
\subsection{Notation}\label{ss2.1}
The ground field is $\C$ and we work in the category of complex
analytic spaces.
A {\em vector bundle} always means a vector bundle of finite rank.
If $E$ is a vector bundle, we denote by $E^*$ the dual vector bundle.
For the tensor power of a line bundle $E$ we write simply $E^k$.
More generally, the convention $E^{-k}=(E^*)^k$ enables us to assume
that $k\in \Z$. By $H^i(E)$ we denote the $i$-th cohomology group
of the sheaf of local sections of $E$, and put $h^i(E)=\dim H^i(E)$,
and  $\chi(E):=\sum_i (-1)^i h^i(E)$.

By $G(k,n)$ we denote the grassmannian of $k$-dimensional subspaces
in $\C^n$. For $k=1$ we have the projective space
$\P^{n-1}:=G(1,n)$. The grassmannian $G(k,n)$ is isomorphic to $G(n-k,n)$ via
the isomorphism
\be{*}
\big\{V\subset \C^n\big\}\mapsto\big\{V^\perp\subset \C^n\big\},
\ee{*}
where $\C^n$ is identified with its dual space.
Under a {\em projective subspace} of $G(k,n)$ we understand
the set of $k$-dimentional subspaces $V\subset\C^n$ such that
$U\subset V\subset W$, where $U\subset W$ are fixed subspaces of $\C^n$
with $\dim U=k-1$, $\dim W>k$, or $\dim W=k+1$, $\dim U<k$. The projective subspace is
a {\it line} if $\dim W - \dim V =2$.

Let $S_k$ be the vector bundle on $G(k,n)$ with  fiber $V$ at the point
$V\in G(k,n)$. There is a canonical inclusion $S_k\subset \tilde{\C}^n$,
where $\tilde{\C}^n$ the trivial bundle on $G(k,n)$ with fiber $\C_n$.
We set $S_{n-k}:=(\tilde{\C}^n/S_k)^*$. By definition, $S_k$ and $S_{n-k}$
are the {\em tautological bundles} on $G(k,n)$. The Picard group
$\Pic(G(k,n))$ is isomorphic to $\Z$, and both maximal exterior powers
$\wedge^kS_k$ and  $\wedge^{n-k}S_{n-k}$ are isomorphic generators
of $\Pic (G(k,n))$. We set $\OO_{G(k,n)}(-1):=\wedge^kS_k\simeq \wedge^{n-k}S_{n-k}$
and $\OO_{G(k,n)}(m):=\OO_{G(k,n)}(-1)^{-m}$ for $m\in\Z$.

An {\em ind-space} is the union $\cup_N X_N$ of analytic spaces $X_N$ related
by closed immersions
\be{}\label{inds}
X_1\tor{\ff_1} X_2\tor{\ff_2}...\ .
\ee{}
For instance, the projective ind-space $\P^\infty$ is the union $\cup_N \P_N$ where $\P^N \subset \P^{N+1}$ is the standard closed immersion.  A {\em morphism of ind-spaces}  $\ff:X\to Y$ is a map whose restriction
to each $X_N$ is a morphism of $X_N$ into $Y_{j(N)}$ for some $j(N)$.
In this paper we will more generally consider sequences (\ref{inds})
of arbitrary morphisms. A {\em vector bundle} on the system (\ref{inds})
is a collection $\{ E_N \}$ of vector bundles $E_N$ on $X_N$ such that $\ff_N^*E_{N+1}=E_N$.
If (\ref{inds}) determines an ind-space $X$, we speak of a vector bundle on $X$.

\subsection{The Barth-Van de Ven-Tyurin Theorem}
\begin{thm}\label{thmGT}
If $X\simeq\P^1$ or $X\simeq\P^\infty$, every vector bundle $E$ on $X$ is
isomorphic to a unique direct sum of line bundles, i.e.
$E\simeq\oplus_j\OO_X(d_j)$ for some unique integers
$d_1\geq...\geq d_{\rk E}$.
\end{thm}
For $X\simeq\P^1$ this is a classical result due to A. Grothendieck.
For $X\simeq\P^\infty$ the theorem has been proved by A. Tyurin in
\cite{T} (and
earlier by W. Barth and A. Van de Ven, \cite{BV}, for vector bundles of rank two). E. Sato also presents a proof in \cite{S1}.

We will make extensive use of Theorem \ref{thmGT}. In particular,
if $l\subset Y$ is a rational curve in a complex manifold or
ind-space $Y$ and $E$ is a vector bundle on $Y$, we call the
numbers $d_1,...,d_{\rk E}$, together with the multiplicities with
which they occur in the isomorphism
$E|_l\simeq\oplus_j\OO_Y(d_j)$, the {\em splitting data} of
$E|_l$. We always denote by $d_1$ the maximal among all $d_j$, and
by $D$ -- the difference $d_1-d_{\rk E}$. When we need to
emphasize the dependence of $d_j$ and $D$ on $l$, we write
$d_j(E|_l)$ and $D(E|_l)$.

We will use the same convention when $Y$ is an ind-space and $l$
is replaced by an ind-subspace $P\subset Y$ isomorphic to $\P^\infty$.
Note that, if $E'$ is a vector bundle on $X\simeq \P^1$
or $X\simeq\P^\infty$ with splitting data $d_1\geq...\geq d_{\rk E}$, and $\nu_1$ is the
multiplicity of $d_1$, there is a unique subbundle $E'_1$ of $E'$
isomorphic to the direct sum of $\nu_1$ copies of $\OO_X(d_1)$.
Indeed, consider the twisted bundle $E'\ot\OO_X(-d_1)$ and
its unique maximal subbundle $E''$ generated by global sections.
Then $E'_1=E''\ot\OO_X(d_1)$.

We say that a vector bundle $E$ on a grassmannian $G(k,n)$ is
{\em uniform}, if the splitting data $d_1(E|_l),...,d_{\rk E}(E|_l)$ for
a line $l\subset G(k,n)$ does not depend on the choice of $l$. The bundle $E$ is {\it linearly trivial} if its restriction $E|_l$ is trivial for any line $l \subset G(k,n)$. A linearly trivial bundle is necessarily trivial: for $k=1$ this is a well-known result, and an induction argument on $k$ yields the result also for $G(k,n)$, see for instance \cite{P}.

\subsection{Local rigidity of direct sums of line bundles on $\P^N$}

Let $\S$ be an analytic space  and $\pi_1$ and $\pi_2$ denote
respectively the projections of $\P^1\times S$ onto the first and
the second factor. If $E$ is a vector bundle on $\P^1\times \S$,
$E(x)$ denotes the restriction of $E$ to $\pi_2^{-1}(x)\simeq
\P^1$. In \cite{D1} the following statement is proved.

\begin{propn}\label{propn1} Let $E$ be a vector bundle on  $\P^1\times \S$,
where $\S$ is a connected  analytic space.
Fix $x_0\in \S$. Then

a) there exists an open neighborhood $U \ni x_0$
such that $d_1(E(x))\leq d_1(E(x_0))$ for $x\in U$, i.e.
$d_1(E(x))$ locally decreases;

b) $\sum_i^{\rk  E(x)} \nu_i(E(x))d_i(E(x))=c_1(E(x))$, where $c_1(E(x))$ is
the first Chern class of $E(x)$ and $\nu_i(E(x))$ is the multiplicity of $d_i(E(x))$
(in particular $\sum_i^{\rk E(x)} \nu_i(E(x))d_i(E(x))$
does not depend on $x\in S$).
\end{propn}
Proposition \ref{propn1}, applied to $E^*$ implies that
$d_{\rk E}(E(x))=-d_1(E^*(x))$ locally increases, hence
$D(E(x))=d_1-d_{\rk E}$ also locally decreases.

\begin{propn}\label{rigid} Let $E$ be a vector bundle on $\P^N$, $N>1$,
which is isomorphic to a direct sum of line bundles. Then $E$ is
locally rigid.
\end{propn}
\begin{proof}
Note that $H^1(\End(E))=0$. This follows from the fact that $\End(E)$ is isomorphic
to a direct sum of line bundles, as $H^1(\LL)=0$ for any line
bundle $\LL$ on $\P^N$ for $N>1$. Therefore the parameter space of the
versal deformation of $E$ consists of a single reduced point, see \cite{D2}.
\end{proof}

\subsection{Estimates related to vector bundles on $\P^1\times\P^1$}
In the case of $\P^1\times \P^1$ we call the fibers of $\pi_1$
(respectively, $\pi_2$) vertical (resp., horizontal)
sections of $\P^1\times \P^1$, and for any vector bundle $E$ on $\P^1\times\P^1$
we denote the twisted bundle $E\ot\pi_1^*(\OO(k))\ot\pi_2^*(\OO(l))$ by
$E(k,l)$.

\begin{lemma}\label{lem1}
Let $E$ be a vector bundle on $\P^1\times \P^1$.
Let $P$ be a vertical section and $\S$ be a horizontal section of
$\P^1\times \P^1$.
Then $d_1(E|_P), d_1(E|_S)\geq -1$ implies
\be{}\label{uneq3}
h^0(E)\leq (\rk E) (d_1(E|_P)+1)(d_1(E|_S)+1).
\ee{}
If $d_1(E|_P)<0$ or $d_1(E|_S)<0$, we have
$h^0(E)=0$.
\end{lemma}

\begin{proof}
The second statement in the Lemma follows immediately from Proposition
\ref{propn1}, since $h^0(F)=0$ for any vector bundle $F$ on $\P^1$
when $d_1(F)<0$.
If $d_1(E|_P)$ or $d_1(E|_S)$
equals $-1$, then (\ref{uneq3}) is obvious, therefore we can assume
that $d_1(E|_P), d_1(E|_S)\geq 0$.

For any $k$ there is the exact sequence
\be{*}
0\tor{} E(0,-k-1)\tor{} E(0,-k)\tor{} E|_S\tor{}0.
\ee{*}
The corresponding cohomology sequence gives
\be{*}\label{uneq1}
h^0(E(0,-k))-h^0(E(0,-k-1))\leq h^0(E|_S)\leq (\rk  E) (d_1(E|_S)+1).
\ee{*}
Summation from $k$ to zero yields
\be{}\label{uneq2}
h^0(E)-h^0(E(0,-k-1))\leq (\rk E) (k+1)(d_1(E|_S)+1).
\ee{}
When $k=d_1(E|_P)$ we have,
$h^0(E(0,-k-1))=0$, therefore (\ref{uneq2}) implies (\ref{uneq3}).
\end{proof}

\begin{lemma}\label{lem2}
Let $E$, $P$, and $S$ be as in Lemma \ref{lem1}.
Suppose $d_1(E|_P), d_1(E|_S)\geq 0$ and
$d_{\rk E}(E|_P), d_{\rk E}(E|_S)\leq 0$.
If $P'$ is another vertical section of $\P^1\times \P^1$,
then
\be{}\label{uneq4}
D(E|_{P'})\leq 4 (\rk E) (D(E|_P)+2)(D(E|_S)+1)-\chi(E)-\chi(E^*).
\ee{}
\end{lemma}

\begin{proof}
From the cohomology sequence of the exact sequence
\be{*}
0\tor{} E\tor{} E(1,0)\tor{} E|_{P'}\tor{}0
\ee{*}
we have
\be{*}
h^0(E)+h^0(E|_{P'})\leq h^0(E(1,0))+h^1(E).
\ee{*}
This, together with the equality $h^2(E)=h^0(E^*(-2,-2))$
(Serre duality), gives
\be{}\label{in1}
d_1(E|_{P'})\leq h^0(E|_{P'})\leq h^0(E(1,0))-\chi(E)+h^0(E^*(-2,-2)).
\ee{}
Note that $d_1((E^*(-2,-2)|_P)=-d_{\rk E}(E|_P)-1$ and
$d_1((E^*(-2,-2)|_S)=-d_{\rk E}(E|_S)-1$, so Lemma \ref{lem1} yields
\be{}\label{in2}
h^0(E^*(-2,-2))\leq (\rk E)(-d_{\rk E}(E|_P)(-d_{\rk E}(E|_S).
\ee{}
Combining (\ref{in1}) and (\ref{in2}) we have
\be{*}
d_1(E|_{P'})\leq (\rk E)(d_1(E|_P)+2)(d_1(E|_S)+1)-\chi(E)
+  (\rk E) (-d_{\rk E}(E|_P))(-d_{\rk E}(E|_S)).
\ee{*}
Since $d_1(E|_P),-d_{\rk E}(E|_P)\leq D(E|_P)$ and
$d_1(E|_S),-d_{\rk E}(E|_S)\leq D(E|_S)$,
we obtain
\be{}\label{uneq4a}
d_1(E|_{P'})\leq 2(\rk E)(D(E|_P)+2)(D(E|_S)+1)-\chi(E).
\ee{}
By adding to (\ref{uneq4a}) the analogous inequality for $E^*$
we obtain (\ref{uneq4}).
\end{proof}

\section{The case of a twisted projective ind-space}
Our first main result is the following theorem.
\begin{thm}\label{thmour}
Let
\be{*}
\P^{i_1}\tor{\ff_{1}}\P^{i_2}\tor{\ff_{2}}...
\ee{*}
be a twisted sequence of projective spaces, i.e.
$\deg{\ff_N}>1$ for infinitely many $N$,
and let $E=\{E_N\}$ be a vector bundle on this inductive system.
Then $E$ is trivial, i.e. all $E_N$ are trivial.
\end{thm}

\begin{lemma}\label{lemChcl}
Let $E$ be as in Theorem \ref{thmour}. Then all
Chern classes $c_q(E_N)$ equal zero.
\end{lemma}

\begin{proof}
By definition, $c_q(E_N)\in H^{2q}(\P^{i_N},\Z)$. For each $j>N$, the homomorphism
\be{}\label{ccl}
H^{2q}(\P^{i_j},\Z)\tor{}H^{2q}(\P^{i_N},\Z)
\ee{}
induced by the composition
$\ff_{j,N} :=  \ff_N\circ...\circ\ff_j$ maps $c_q(E_j)$ to $c_q(E_N)$.
On the other hand, we can identify $H^{2q}(\P^{i_N},\Z)$ and
$H^{2q}(\P^{i_j},\Z)$ in a standard way with $\Z$, and (\ref{ccl})
considered as an endomorphism of $\Z$ is nothing but multiplication
by $(\deg \ff_j)\cdot...\cdot(\deg\ff_N)$. Therefore $c_q(E_N)\in\Z$ is divisible
by $(\deg\ff_j)\cdot...\cdot(\deg\ff_N)$ for all $j>N$.
According to our assumption, $\lim_{j\to\infty}((\deg\ff_j)\cdot...\cdot(\deg\ff_N))=\infty$,
consequently, $c_q(E_N)=0$.
\end{proof}

\noindent{\em Proof of Theorem \ref{thmour}}.
Without restriction of generality we assume that $i_1=1$.
Fix $N$ with $i_N>2$ and set
\be{*}
D_N:=\max_{l'}D(E_N|_{l'}),
\ee{*}
where $l'$ runs over all  lines in $\P^{i_N}$. Let $l$ be a  line in $\P^{i_N}$
such that $D(E|_l)=D_N$, and $Q$  be a projective subspace in $\P^{i_N}$
of dimension $i_{N}-2$ not intersecting $l$.
We can choose homogeneous
coordinates $z_0,...,z_{i_N}$ in $\P^{i_N}$ such that $l$
is defined by $z_i=0$, $i\geq 2$,
and $Q$ is defined by $z_0=z_1=0$.
Furthermore, we fix a morphism $f:\P^1\to \P^{i_N}$ of degree $\deg \ff_{1,n} $ such that
\be{}\label{intersec}
f(\P^1)\cap l=\emptyset, \quad f(\P^1)\cap Q=\emptyset,
\ee{}
and
\be{}\label{uneqd1}
D(f^*E_N)\leq D(E_1).
\ee{}
To see that $f$ exists, note that by Proposition \ref{propn1}
there is a neighborhood of unity $U$ in the group
of linear transformation of $\P^{i_N}$ such that for any $g\in U$
the morphism $g\circ \ff_{1,N}:\P^1\to\P^{i_N}$ (where $\ff_{1,N}=\ff_{N}\circ \ff_{N-1}\circ ... \circ \ff_1$)
 satisfies (\ref{uneqd1}).
Moreover, it is obvious that one can choose $g$
in such a way that (\ref{intersec}) holds.

We now extend the morphism $f$ to a morphism
$\tilde f:\P^1\times \P^1\to \P^{i_N}$ in the following way.
Let $f$ be given in coordinates as
$z_i=f_i(x_0,x_1)$, where $x_0, x_1$ are homogeneous
coordinates on the space $\P^1$ which we identify with the second factor
of $\P^1\times \P^1$. Let $t_0, t_1$ be  homogeneous
coordinates on the first factor of $\P^1\times \P^1$.
We define $\tilde f$ by putting
$z_i=t_0f_i(x_0,x_1)$ for $i=0,1$ and
$z_i=t_1f_i(x_0,x_1)$ for $i\geq 2$.
Condition (\ref{intersec}) ensures that
$\tilde f$ is well-defined.
Set $\tE_N:=\tilde f^*_N E_N$.

The morphism $\tilde f$ has the following properties.

a) The restriction of $\tilde f$ to the vertical fiber $P$ over the point
$(t_0,t_1)=(1,1)$ coincides with $f$.
Therefore (\ref{uneqd1}) implies
\be{}\label{uneqP}
D(\tE_N|_P)&\leq & D({E_1}).
\ee{}
b) The restriction of $\tilde f$ to the vertical fiber $P'$
over the point $(t_0,t_1)=(1,0)$ is a morphism
of $P'$ onto $l$ of degree $\deg \ff_{1,N}$.
Thus
\be{}\label{eqP}
D(\tE_N|_{P'})&= & (\deg\ff_{1,N})D_N.
\ee{}
c) The restriction of $\tilde f$ to any horizontal fiber $S$
is a linear embedding. Therefore
\be{}\label{uneqS}
D(\tE_N|_S)&\leq & D_N.
\ee{}

The hypotheses of Lemma \ref{lem2} hold for $\tE_N$, since
$c_1(\tE_N|_{P})=c_1(\tE_N|_{S})=0$ by Lemma \ref{lemChcl}. Therefore
Lemma \ref{lem2}, together with (\ref{uneqP}),
(\ref{eqP}), and (\ref{uneqS}), yields
\be{}\label{muneq}
(\deg \ff_{1,N})D_N\leq 4 (\rk E_N)(D({E_1})+2)(D_N+1)-\chi(E_N)-\chi(E^*_N).
\ee{}
By the Riemann-Roch Theorem and the vanishing of all Chern
classes of $E_N$ (Lemma \ref{lemChcl}), we have
\be{*}
\chi(E_N)=\chi(E^*_N)=n\chi(\O_{P^N})=\rk E_N.
\ee{*}
Therefore (\ref{muneq}) turns into
\be{}\label{cuneq}
(\deg\ff_{1,N})D_N\leq 4 (\rk E_N)(D(E_1)+2)(D_N+1)-2 \rk E_N.
\ee{}
As $\rk E_N = \rk E$ and $D(E_1)$ do not depend on $N$, we see that the right hand side of (\ref{cuneq})
is a linear function of $D_N$, while the left hand side
grows with $N$ faster than a linear function of $D_N$, as $\lim_{N\to\infty}\deg \ff_{1,N}=\infty$ by our hypothesis.
Therefore, inequality (\ref{cuneq}) can hold
only if $D_N$ equals zero for large enough $N$ and hence
for all $N$. This means that all vector bundles $E_N$ are linearly
trivial, and thus trivial. \qed

Theorem \ref{thmour} together with Theorem \ref{thmGT} leads
to the following general description of a vector bundle on an
inductive limit of relative projective spaces. Let $p_N:M_N\to S_N$ be
a relative projective space, i.e. a locally trivial fibration with
base a connected complex analytic space $S_N$ and with fiber the
projective space $\P^{i_N}$. Let furthermore
\be{}\label{reld}
\begin{CD}
M_1@>\ff_1>>M_2@>\ff_2>>...\\
@Vp_1VV @Vp_{2}VV\\
S_1@>>>S_2@>>>...
\end{CD}
\ee{}
be a commutative diagram.
Note that if the restriction of $\ff_N$ to a fiber $P_N$ of $p_N$ has degree $d$, then
the restriction of $\ff_N$ to any fiber of $p_N$ has also degree $d$.
Therefore it makes sense to call (\ref{reld}) {\em linear} if the degree of the restriction
of $\ff_N$ on the fibers of $p_N$ equals one for almost all $N$,
and {\em twisted} otherwise.
\begin{thm}\label{relthm1}
Let $E=\{E_N\}$ be a vector bundle on the upper row of (\ref{reld}).
Then there exist integers $d_1\geq...\geq d_{\rk E}$ such that,
for large enough $N$ and any fiber $P_N$ of $p_N$, $E_N|_{P_N}\simeq\oplus_r\OO_{P_N}(d_r)$.
Furthermore, if (\ref{reld}) is twisted, $d_1=...= d_{\rk E}=0$.
\end{thm}

\begin{proof}
The fact that $E_N|_{P_N}$ is isomorphic to a direct sum of line bundles follows
from Theorems \ref{thmGT} and \ref{thmour}.
The fact that the splitting data of $E_N|_{P_N}$ does not depend on $P_N$
follows from Proposition \ref{rigid}.
If (\ref{reld}) is linear, the degree of $P_N$ equals one for  large enough $N$,
thus the splitting data $E_N|_{P_N}$ does also not depend on $N$ for large
enough $N$. Finally, if (\ref{reld}) is twisted, each $d_r(E_N|_{P_N})$ is
divisible by $(\deg\ff_N)\cdot...\cdot(\deg\ff_j)$ for all $j>N$,
and therefore $d_r(E_N|_{P_N})=0$.
\end{proof}
In general it is not true that $E$ is isomorphic to a direct sum of line bundles, or
that $E$ admits a filtration whose associated quotient are line bundles.
The reader will easily prove that the latter holds only under
the additional assumption that all $d_r$ are distinct.
Finally, if (\ref{reld}) is twisted, Theorem \ref{relthm1}
implies that each $E_N$ is the pull-back
$p^*_NE'_N$ for some vector bundle $E'=\{E'_N\}$ on the lower row of
(\ref{reld}).

\section{The case of ind-grassmannians}

In this section we consider two different types of closed immersions
\be{}\label{seqgr}
G(k_1,n_1)\tor{f_1} G(k_2,n_2)\tor{f_2}...
\ee{}
and characterize vector bundles on the corresponding ind-spaces.

\subsection{Standard extensions of grassmannians}\label{ss4.1}

We define a {\em standard extension} of grassmannians as a closed immersion of the form

\be{*}
\lambda_{r,m}:G(k,n)\to G(k+r,n+m), \qquad
\{V\subset\C^n\}\mapsto\{V\oplus W\subset\C^n\oplus\C^m\},
\ee{*}
where $W \subset \C^m$ is a fixed subspace of dimension $r\geq 0$.

\begin{propn}\label{proplevel}
Assume that a sequence (\ref{seqgr}) is given, where $f_N$ are standard extensions.
Let $E=\{E_N\}$ be a vector bundle on (\ref{seqgr}).
Then for each $N$, $E_N$ is a uniform
bundle on $G(k_N,n_N)$.
\end{propn}
\begin{proof}
For each $N$ consider the natural  diagram
\be{}\label{Di}
\begin{CD}
F(k_N-1,k_N,n_N) @>\pi_N>> G(k_N,n_N)\\
@Vp_NVV \\G(k_N-1,n_N),
\end{CD}
\ee{}
where $F(k_N-1,k_N,n_N)$ stands for the space of all flags of type ($k_N-1, k_N$) in  $\C^{n_N}$. Note that $f_N$ induce morphisms between the diagrams (\ref{Di})
for $N$ and $N+1$. Furthermore, the embeddings
$F(k_N-1,k_N,n_N)\to F(k_{N+1}-1,k_{N+1},n_{N+1})$ and the morphisms $p_N$
define a relative projective ind-space. The vector bundles
$\pi^*_N E_N$ define a vector bundle on this relative projective ind-space.
By Theorem \ref{relthm1} the splitting data of the restriction
$\pi^*_N E_N|_{P_N}$ on each fiber $P_N$ of $p_N$ does not depend on $N$ and $P_N$.
This implies the result, as any projective line $l\subset G(k_N,n_N)$ is a line
in some fiber $P_N$, and the splitting data of $E|_l$ and $\pi^*_NE_N|_{P_N}$
are equal.
\end{proof}

\begin{thm}\label{thmgr1}
Assume that a sequence (\ref{seqgr}) is given,
where $f_N$ are standard extensions and
$\lim_{N \to \infty} k_N = \lim_{N \to \infty} (n_N - k_N) = \infty$.
Then any vector bundle $E=\{E_N\}$ on the inductive limit of (\ref{seqgr})
is isomorphic to  a direct sum of line bundles.
\end{thm}
\begin{proof}
By Proposition \ref{proplevel}, each $E_N$ is a uniform bundle.
We will prove the Theorem by induction on $\rk E = \rk E_N$.
Fix $N$ and let $d_1,...,d_{\rk E_N}$ be the splitting  data
for the restriction $E_N|_l$, where $l  \subset G(k_N, n_N)$ is any line.
Denote by $\nu_1$ the multiplicity of $d_1$. As each $f_j$ induces an
isomorphisms of Picard groups, the splitting data $d_1,...,d_{\rk E_N}$
does not depend on $N$. Fix $x \in G(k_N, n_N)$.
Any  line $l \subset G(k_N, n_N)$ passing through $x$ determines a
subspace $E_{l,x}'$ in the fiber $(E_N)_x$: this is the fiber of the unique
subbundle  of $E_N|_l$  isomorphic to the direct sum of $\nu_1$ copies
of $\OO_l(d_1)$, see Subsection 2.2. The  variety of all lines $l$ passing
through $x$ is isomorphic to the direct product
$\P^{k_N  -1}  \times \P^{n_N - k_N -1}$, hence we have a morphism
\be{*}
\psi_x : \P^{k_N  -1}  \times \P^{n_N - k_N -1} \tor{} G(\nu_1, \rk E_N), \qquad
l \mapsto E_{l,x}'.
\ee{*}
By our assumption, $\lim_{N \to \infty} k_N = \lim_{N \to \infty} (n_N - k_N) = \infty$.
Therefore, for large $N$, the morphism $\psi_x$ is trivial,
i.e. it determines a fixed subspace $(E'_N)_x$ in  $(E_N)_x$.
In this way we obtain a subbundle $E_N' \subset E_N$.
Clearly $f_N^*(E'_{N+1}) = E'_N$, therefore $E' := \{ E'_N\}$
is a well-defined subbundle of $E$.

If $\rk E' = \rk E$, twisting by the line bundle $\{\OO_{ G({k_N},n_N)}(-d_1)  \}$ yields
a bundle whose restriction to any line in  $G(k_N, n_N)$ is trivial
for all $N$. Therefore the twisted bundle is trivial on each $G(k_N, n_N)$,
and hence trivial. Thus $E$ is isomorphic to a direct sum of $\nu_1$
copies of $\{\OO_{ G({k_N},n_N)}(d_1)  \} $. If $\rk E' < \rk E$,
the induction assumption implies that both $E'$ and $E/E'$ are isomorphic
to direct sums of line bundles. Finally, the observation that there are
no non-trivial extensions of line bundles on $G(k_N,  n_N)$ for $k_N>1$
(by the Bott-Borel-Weil Theorem, $H^1( \LL ) =0$ for
any line bundle $\LL$ on $G(k_N, n_N)$ unless $k_N = n_N - k_N = 1$)
implies that $E$ is isomorphic to a direct sum of line bundles.
\end{proof}

Here is an important special case of  Theorem \ref{thmgr1}.
If $V\subset\C^\infty$ is an arbitrary fixed subspace, a subspace
$V'\subset\C^\infty$ is {\it commensurable} with $V$, if there exists a finite
dimensional subspace $U\subset\C^\infty$ such that $V\subset V'+U$,
$V'\subset V+U$, and $\dim V\cap U=\dim V'\cap U$.
The set of all subspaces $V'$ comensurable with $V$ is, by definition,
the ind-grassmannian $G(V,\infty)$. If $V$ is finite dimensional and
$\dim V=k$, then $G(V,\infty)=G(k,\infty)$. In \cite{DP}
an explicit construction of $G(V,\infty)$ as an ind-space is given.
Moreover, $G(V,\infty)$ is the inductive limit of standard extensions,
and the conditions of Theorem \ref{thmgr1} are satisfied for $G(V,\infty)$
if and only if $\dim V=\codim_{\C^\infty}V=\infty$. Therefore,
by Theorem \ref{thmgr1}, we conclude that in the latter case
every vector bundle on $G(k,\infty)$ is isomorphic to a direct sum of line bundles.

More generally, if $V_1\subset...\subset V_r$ is any fixed flag in $\C^\infty$,
the ind-variety $F(V_1,...,V_r,\infty)$ of flags $V'_1\subset...\subset V'_r$
comensurable with $V_1\subset...\subset V_r$, is constructed in \cite{DP}.
We leave it to the reader to prove the following corollary of Theorem \ref{thmgr1}
by double induction on $\rk E$ and $r$.
\begin{cor}
Let $\dim V_1=\dim (V_2/V_1)=...=\dim(\C^\infty/V_r)=\infty$.
Then every vector bundle $E$ on $F(V_1,...,V_r,\infty)$ is isomorphic
to a direct sum of line bundles.
\end{cor}

\subsection{Twisted extensions of grassmannians}

An alternative definition of a grassmannian is as a homogeneous space $GL(n)/P$ for a maximal parabolic
subgroup $P \subset GL(n)$. We call a morphism
$f:G(k_1,n_1)\to G(k_2,n_2)$ {\em  homogeneous}
if it is induced by a group
homomorphism $\tilde f:GL(n_1)\to GL(n_2)$. A homogeneous morphism is a closed
immersion or its image is a point.

Consider a homogeneous morphism $f : G(k_1, n_1) \to G(k_2, n_2)$, where $G(k_1,n_1)=GL(n_1)/P_1$,
$G(k_2,n_2)=GL(n_2)/P_2$, and $P_1$, $P_2$ are maximal parabolic subgroups
such what $\tilde f(P_1)\subset P_2$. The reductive part of $P_1$ is isomorphic to
$GL(k_1)\times GL(n_1-k_1)$ and the reductive part of $P_2$ is isomorphic to
$GL(k_2)\times GL(n_2-k_2)$. Note that $\tilde f(GL(k_1)\times GL(n_1-k_1))$ is
contained in one component of the direct product $GL(k_2)\times GL(n_2-k_2)$
if and only if the induced morphism $f:G(k_1,n_1)\to G(k_2,n_2)$ is a morphism
into a point. In the sequel we adopt a convention: whenever we write a homogeneous
morphism as $f:G(k_1,n_1)\to G(k_2,n_2)$ and $f$ is a closed immersion,
we automatically assume that $\tilde f(GL(k_1))\subset GL(k_2)$ and
$\tilde f(GL(n_1-k_1))\subset GL(n_2-k_2)$. Furthermore, we call a morphism
$\ff : \P^N \to G(k,n)$ {\em $k$-split} if
$\ff^*S_k$ is isomorphic to a direct sum
of line bundles on $\P^N$.
For example, Theorem \ref{thmGT} implies that any morphism $\P^1\to G(k,n)$ is $k$-split.

\begin{propn}\label{gmap}
a) Any homogeneous morphism  $\ff:\P^N = G(1, N) \to G(k,n)$ is $k$-split.

b) Let $f:G(k_1,n_1)\to G(k_2,n_2)$ be a homogeneous embedding of grassmannians
and $\ff:\P^N\to G(k_1,n_1)$ be a $k_1$-split morphism. Then the
composition $f \circ \ff:\P^N\to G(k_2,n_2)$ is $k_2$-split.
\end{propn}
\begin{proof}
a) Let $\P^N=GL(N+1)/P_1$, and $G(k,N) = GL(n)/P_2$. The reductive parts of $P_1$
and $P_2$ are isomorphic respectively
to $GL(1) \times GL(N)$ and $GL(k)\times GL(n-k)$. According to our convention,
the group homomorphism $\tilde{\ff}:GL(N+1) \to GL(n)$ which induces $\ff$ maps $GL(1)$
into $GL(k)$, and  $GL(N)$ into $GL(n-k)$.
Therefore the structure group of $\ff^* S_k$ is reduced to $GL(1)$,
which implies that $\ff^* S_k$ is a direct sum of line bundles.

b) Let $\tilde{f} : GL(n_1) \to GL(n_2)$ be a group homomorphism which
induces the morphism $f$. Since $\tilde{f} (GL(k_1)) \subset GL(k_2)$, the
$GL(n_1)$-homogeneous bundle $f^*S_{k_2}$ is a direct summand in the tensor
product of several copies of $S_{k_1}$ and $S^*_{k_1}$. This,
together with the fact, that the morphism $\ff$ is $k_1$-split,
implies that $(f \circ \ff)^*S_{k_2}$ is a direct summand of a
bundle isomorphic to a direct sum of line bundles on $\P^N$.
Hence  $(f \circ \ff)^*S_{k_2}$ is isomorphic itself to a
direct sum of line bundles, i.e. the morphism $f \circ \ff$ is $k_2$-split.
\end{proof}

Here is a coordinate form of a $k$-split map.
Let $t=(t_0,...,t_N)$ be homogeneous coordinates on $\P^N$,
and $v_1,...,v_n$ be a basis in $\C^n$. Then the reader will check that  a morphism
$\ff:\P^N\to G(k,n)$ is $k$-split
if and only if $\ff$ can be
presented in the form
\be{}\label{span}
\ff(t)=\span\Big\{\sum_{j=1}^n\ff_{i,j}(t)v_j, \ i=1,...,k\Big\}
\ee{}
for some homogeneous polynomials   $\ff_{i,j}(t)$  in $t$.
It is clear that if   $\ff:\P^N\to G(k,n)$ is represented in the form
(\ref{span}), then $\deg\ff_{i,j}$ does not depend on $j$ and
$\deg\ff=\sum_{i=1}^k\deg\ff_{i,j}$.

\begin{lemma}\label{lem11}
Let $\ff:\P^N\to G(k,n)$ be a $k$-split morphism and
$\lambda_{0,k}:G(k,n)\to G(k,n+k)$ be a standard extension
(see Subsection \ref{ss4.1}).
Then the composition $\lambda_{0,k}\circ\ff:\P^N\to G(k,n+k)$
factors through a standard extension $\P^N\to\P^{N+1}$, via a $k$-split morphism $\ff:\P^{N+1} \to G(k,n+k)$.
\end{lemma}
\begin{proof}
Let $\ff$  be presented in the form (\ref{span}), and let the
vectors $v_{n+1},...v_{n+k}$ extend the set  $v_1,...,v_n$ to a basis
in $\C^{n+k}$. Fix homogeneous coordinates $t'=(t_0,...,t_{N+1})$ in $\P^{N+1}$ and identify $\P^N$ with the hyperplane  $t_{N+1}=0$. Then $\tilde{\ff} : \P^{N+1} \to G(k, n+k)$ is given by the formula
$$\tilde{\ff}(t')=\span\Big\{\sum_{j=1}^n\ff_{i,j}(t)v_j+t^{\delta_i}_{N+1}v_{n+i},\ i=1,...,k\Big\},$$
where $\delta_i := \deg \ff_{i,j}$.
\end{proof}

We need one last definition.
We call a homogeneous morphism $f:G(k_1,n_1)\to G(k_2,n_2)$
a {\em twisted extension}
if $f$ decomposes as
$\lambda_{0,n_2-n}\circ h$, where
$h:G(k_1,n_1)\to G(k_2,n)$
is a homogeneous morphism for $n_2-n\geq k_2$, and
$\lambda_{0,n_2-n}:G(k_2,n)\to G(k_2,n_2)$ is a standard extension.

\begin{propn}\label{prop100}
Assume that a sequence (\ref{seqgr}) is given, where $f_N$ are twisted extensions,
infinitely many of which have degree greater than one.
Suppose in addition that  $\lim_{N\to\infty}k_N=\lim_{N\to\infty}(n_N-k_N)=\infty$.
Then any vector bundle $E = \{ E_N\} $ on the inductive limit of (\ref{seqgr}) is trivial.
\end{propn}

\begin{proof} Let $E=\{E_N\}$ be a vector bundle on
the inductive limit of (\ref{seqgr}).
Fix $N$ and a $(k_N - 1)$-dimensional subspace of $\C^{n_N}$.
This determines a projective subspace $\P^{n_N - k_N}$ of $G(k_N,n_N)$,
and the closed immersion
$\ff_N : \P^{n_N - k_N} \to G(k_N, n_N)$ is
$k_N$-split as $\ff^* S_{k_N} \simeq (k_N -1)\OO_{\P^{n_N - k_N}} \oplus \OO_{\P^{n_N - k_N}}(-1) $. Furthermore, let $f_N = \lambda_{0,m_N} \circ h_N$, where $h_N: G(k_N,n_N)\to G(k_{N+1},n_{N+1}-m_N)$ is a  homogeneous closed immersion and $m_N \geq k_{N+1}$. Proposition \ref{gmap} b) implies that the composition $h_N \circ \ff_N$ is $k_{N+1}$-split. Then, by Lemma \ref{lem11}, the morphism $f_N \circ \ff_N = \lambda_{0,m_N} \circ (h_N \circ \ff_N )$ factors through a standard extension $\P^{n_N - k_N} \to \P^{n_N - k_N +1}$ via a homogeneous immersion  $\ff_{N+1} : \P^{n_N - k_N+1} \to G(k_{n+1}, n_{N+1})$.  Proceeding by induction on $j$,
we build homogeneous immersions  $\ff_{N+j}:\P^{n_N-k_N+j}\to G(k_{N+j}, n_{N+j})$
and standard extensions $\P^{n_N-k_N+j} \to \P^{n_N-k_N+j+1}$ which form the commutative diagram
\be{}\label{cov}
\begin{CD}
\P^{n_N-k_N}@>>>\P^{n_N-k_N+1}@>>>...\\
@V\ff_N VV @V\ff_{N+1}VV\\
G(k_N, n_N)@>f_N>>G(k_{N+1}, n_{N+1})@>f_{N+1}>>... \ .
\end{CD}
\ee{}
Applying the BVT theorem to the inductive limit of the upper row of (\ref{cov}), we conclude that $E_N|_{\P^{n_N-k_N}}$ is isomorphic to a direct sum of line bundles. Furthermore, one sees exactly as in Proposition \ref{proplevel} that $E_N$ is a uniform bundle. Now the argument in the proof of  Theorem \ref{thmgr1} shows that, since  $\lim_{N\to\infty}k_N=\lim_{N\to\infty}(n_N-k_N)=\infty$,
$E_N$ is isomorphic to a direct sum of line bundles. Therefore $E$ is also isomorphic to a direct sum of line bundles. Finally, as $\deg f_N >1$ for
infinitely many $N$, a line bundle on this inductive limit is  necessarily trivial, hence $E$ is trivial.
\end{proof}

Here are two examples to Proposition \ref{prop100}.

Let $\C^w$ be a proper subspace of $\C^m$. Consider the closed immersion
\be{}\label{ex100}
G(k,n)\tor{}G(kw,nm), \qquad
\{V\subset\C^n\}\mapsto\{V\ot \C^w\subset\C^n\ot\C^m\}.
\ee{}
It factors, via a standard extension, through the closed immersion
\be{*}\label{ex100w}
G(k,n)\tor{}G(kw,nw), \qquad
\{V\subset\C^n\}\mapsto\{V\ot \C^w\subset\C^n\ot\C^w\}.
\ee{*}
Fixing a sequence of pairs $\C^{w_N}\subset\C^{m_N}$, determines a sequence of
immersions (\ref{ex100}). If the set $\{w_N\}$ is
bounded, for sufficiently large $N$ all embeddings in this sequence are
twisted extensions as
$$nm_1...m_{N+1}-nw_1...w_{N+1}\geq kw_1...w_{N+1}$$
for large enough $N$.
Therefore Proposition \ref{prop100} implies that any vector bundle
on the inductive limit is trivial.

Finally consider a closed immersion of the form
\be{}\label{ex101}
G(k,n)\tor{}G(k(k+1)/2,n^2), \qquad
\{V\subset\C^n\}\mapsto\{S^2(V)\subset\C^n\ot\C^n\},
\ee{}
where $S^2$ stands for symmetric square. By iteration of (\ref{ex101}) we obtain
another example of a sequence
of twisted extensions. Indeed,
(\ref{ex101}) factors, via a standard extension, through the closed
immersion
\be{*}
G(k,n)\tor{}G(k(k+1)/2,n(n+1/2), \qquad
\{V\subset\C^n\}\mapsto\{S^2(V)\subset S^2(\C^n)\},
\ee{*}
and the inequality
\be{*}
k(k+1)/2\leq n^2-n(n+1)/2
\ee{*}
 needed in the definition of a twisted extension, is trivially satisfied.
Proposition \ref{prop100} implies that any vector bundle
on the corresponding inductive limit is also trivial.

\section{Conclusion}

The results of this paper lead naturally to the following question.
Consider an arbitrary sequence of morphisms
\be{}\label{con1}
G(k_1,n_1)\tor{f_1} G(k_2,n_2)\tor{f_2}...\ .
\ee{}
Is it true that, if $E= \{ E_N\} $ is an arbitrary vector bundle on (\ref{con1})
then each $E_N$ is a homogeneous bundle on $G(k_N,n_N)$?
In fact, the BVT Theorem, both main theorems of Sato (see \cite{S2}),
Theorems \ref{thmour} and \ref{thmgr1}, as well as
Proposition \ref{prop100}, are all equivalent to
the affirmative answer to the above question for the  cases they apply to.
More precisely, the description of vector bundles on sequences (\ref{con1})
in each specific case considered in the above statements, is equivalent
to the descriptions of systems of homogeneous bundles $\{E_N\}$
with $f_{N+1}^*E_{N+1}=E_N$. Note that Theorem \ref{thmour}
applies to not necessarily homogeneous morphisms, while in
all other statements  homogeneous morphisms are considered.
It would be very interesting to give an answer to the above
general question even under the assumption that all morphisms $f_N$
are homogeneous.

J. Donin\\
Department of Mathematics\\
Bar-Ilan University\\
52900 Ramat Gan, Israel\\
donin@macs.biu.ac.il

\medskip

I. Penkov\\
Department of Mathematics\\
University of California Riverside\\
Riverside, CA 92521, USA\\
penkov@math.ucr.edu

\end{document}